\documentclass[11pt]{article}

\usepackage{amsmath}
\usepackage{amsfonts}
\usepackage{amssymb}
\usepackage{amsthm}
\usepackage{enumitem}
\usepackage{setspace}
\usepackage{etoolbox}

\newtheoremstyle{mystyle}
{11pt}                          
{11pt}                          
{}                                      
{}                                      
{\bfseries}                     
{.}                                     
{5.5pt}                         
{}                                      

\newtheoremstyle{mystyle2}
{11pt}                          
{11pt}                          
{}                                      
{}                                      
{\bfseries}                     
{}                                     
{5.5pt}                         
{}                                      

\theoremstyle{mystyle}
\newtheorem{theorem}{Theorem}

\newtheorem{proposition}{Proposition}
\newtheorem{corollary}{Corollary}

\theoremstyle{mystyle2}

\setitemize{noitemsep, topsep=5.5pt, parsep=5.5pt, partopsep=0pt}

\allowdisplaybreaks
\appto\normalsize{
        \abovedisplayskip=5.5pt plus 2pt minus 2pt
        \belowdisplayskip=5.5pt plus 2pt minus 2pt
        \abovedisplayshortskip=5.5pt plus 2pt minus 2pt
        \belowdisplayshortskip=5.5pt plus 2pt minus 2pt}
\appto\small{
        \abovedisplayskip=5.5pt plus 2pt minus 2pt
        \belowdisplayskip=5.5pt plus 2pt minus 2pt
        \abovedisplayshortskip=5.5pt plus 2pt minus 2pt
        \belowdisplayshortskip=5.5pt plus 2pt minus 2pt}

\newcommand{\gap}{\vspace{11pt}}

\newcommand{\tr}{\operatorname{tr}}

\newcommand{\lchi}{\raisebox{2pt}{$\chi$}}
\newcommand{\R}{{\cal R}}
\newcommand{\Rn}{{\cal R}^n}

\newcommand{\Sn}{{\cal S}^n}

\newcommand{\Hn}{{\cal H}^n}

\newcommand{\wl}{\widetilde{\lambda}}
\newcommand{\G}{{\cal G}}
\newcommand{\V}{{\cal V}}
\newcommand{\W}{{\cal W}}
\newcommand{\lc}{\lambda(c)}

\newcommand{\lx}{\lambda(x)}
\newcommand{\lu}{\lambda(u)}
\newcommand{\ly}{\lambda(y)}

\title{\bf
Commutation principles for optimization problems on spectral sets 
in   Euclidean Jordan algebras}
\author{
        M. Seetharama Gowda\\
        Department of Mathematics and Statistics\\
        University of Maryland, Baltimore County\\
        Baltimore, Maryland 21250, USA\\
        gowda@umbc.edu
}

\date{\today }

\begin{document}

\maketitle

\begin{abstract}
The commutation principle of Ram{\'i}rez, Seeger, and Sossa \cite{ramirez-seeger-sossa} proved in the setting of
Euclidean Jordan algebras says that when the sum of a real valued function $h$ and a
spectral function $\Phi$ is minimized/maximized  over a spectral set $E$,  any local optimizer $a$ at which $h$ is  Fr\'{e}chet differentiable 
operator commutes with the
derivative $h^{\prime}(a)$.
In this  paper, assuming the existence of a subgradient in place  the derivative (of $h$), we establish `strong operator commutativity' relations:
If $a$ solves the problem $\underset{E}{\max}\,(h+\Phi)$, then $a$ strongly operator commutes with every element in the subdifferential of $h$ at $a$; If $E$  and $h$ are convex
and $a$ solves the problem $\underset{E}{\min}\,h$, then $a$ strongly operator commutes with the negative of some element in the subdifferential of $h$ at $a$.
These results improve known (operator) commutativity  relations for linear $h$ and for solutions of variational inequality problems.
We establish these results via  a geometric commutation principle that is valid not only in Euclidean Jordan algebras, but also in the broader setting of FTvN-systems.
\end{abstract}

\vspace{1cm}
\noindent{\bf Key Words:} Euclidean Jordan algebra,   spectral sets/functions,
 commutation principle,  variational inequality problem, normal cone, subdifferential.

\gap

\noindent{\bf MSC2020 subject classification:}   17C20, 17C30, 52A41, 90C26.
\newpage

\section{Introduction}
      {\it Let $\V$ be a Euclidean Jordan algebra of rank $n$ carrying the trace inner product} \cite{faraut-koranyi} 
and $\lambda:\V\rightarrow \Rn$ denote the eigenvalue map (which takes $x$ to $\lambda(x)$, the vector of eigenvalues of $x$ with entries written in the decreasing order). For any $a\in \V$, we define its $\lambda$-orbit by
$[a]:=\{x\in \V: \lambda(x)=\lambda(a)\}.$
A  set $E$ in $\V$ is said to be a {\it spectral set}  if it is of the form
$E = \lambda^{-1}(Q) $
for some (permutation invariant) set $Q$ in $\Rn$ or, equivalently, a union of $\lambda$-orbits.
A function $\Phi:\V\rightarrow \R$ is said to be a {\it spectral function} if it is of the form
$\Phi = \phi \circ \lambda $
for some (permutation invariant)  function $\phi : \Rn \rightarrow \R$. 

\gap

In \cite{ramirez-seeger-sossa},  Ram{\'i}rez, Seeger, and Sossa prove the following commutation principle in $\V$: 
\begin{theorem}\label{ramirez theorem} 
{\it Suppose $a$ is a local optimizer of the problem
$\underset{E}{\min\hspace{-1mm}/\hspace{-1mm}\max}\,( h+\Phi),$
where $E$ is a spectral set, $\Phi$ is a spectral function, and $h:\V\rightarrow \R$. If $h$ is Fr\'{e}chet differentiable at $a$, then $a$ and $h^{\prime}(a)$ operator commute.}
\end{theorem}

In \cite{gowda-jeong}, Gowda and Jeong extended the above result 
by assuming that  $E$ and $\Phi$ are invariant under the automorphisms of $\V$  and stated an 
analogous result in the setting of normal decomposition systems.
Subsequently, certain  modifications (such as replacing the sum by other combinations) and applications were given 
by Niezgoda \cite{niezgoda}.

\gap

The main objective of this note is to describe some analogs of the above commutation principle by assuming the existence of a subgradient in place of the derivative  (of $h$).
In each analog, this change results in a  stronger commutativity relation. We derive
 these analogs  via a geometric commutation principle. 
To elaborate,  we first recall some definitions.  

\begin{itemize}
\item [$\bullet$] 
 We say that elements $a$ and $b$ {\it operator commute} in $\V$
if there exists a Jordan frame $\{e_1,e_2,\ldots, e_n\}$ in $\V$ such that the spectral decompositions of $a$ and $b$ are given by 
$$a=a_1e_1+a_2e_2+\cdots+a_ne_n\quad \mbox{and}\quad b=b_1e_1+b_2e_2+\cdots+b_ne_n,$$
where $a_1,a_2,\ldots, a_n$ are the eigenvalues of $a$ and $b_1,b_2,\ldots, b_n$ are the eigenvalues of $b$. 
If, additionally, above decompositions hold with $a_1\geq a_2\geq \cdots\geq a_n$ and $b_1\geq b_2\geq \cdots\geq b_n$, we say that $a$ and $b$ {\it strongly operator commute} (also said to `simultaneously diagonalizable' \cite{lim et al} or said to have `similar joint decomposition' \cite{baes}). Some equivalent formulations are described in the next section.
\item [$\bullet$]
Given a (nonempty) set $S$ in $\V$ and $a\in S$, we define the {\it normal cone} of $S$ at $a$ by 
$$N_{S}(a):=\{d\in \V: \langle d,x-a\rangle \leq 0\,\,\mbox{for all}\,\,x\in S\}.$$
\item [$\bullet$] Let  $h:\V\rightarrow \R\cup \{\infty\}$, $S\subseteq \V$, and $a\in S\cap \mbox{dom}\,h$.  We define the {\it subdifferential of $h$ at $a$ relative to $S$} by 
$$\partial_{S}\,h(a):=\{d\in \V: h(x)-h(a)\geq \langle d,x-a\rangle\,\,\mbox{for all}\,\,x\in S\};$$ any element of $\partial_{S}\,h(a)$ will be called
 a {\it $S$-subgradient of $h$ at $a$}.
Finally, when $S=\V$, we define the {\it subdifferential of $h$ at $a$} by 
$$\partial\,h(a):=\{d\in \V: h(x)-h(a)\geq \langle d,x-a\rangle\,\,\mbox{for all}\,\,x\in \V\}.$$ 
We note that subdifferentials  may be empty and $\partial\,h(a)\subseteq \partial_{S}\,h(a)$.
We also note \cite{rockafellar} that when $h$ (defined on all of $\V$) is convex, the subdifferential is  nonempty, compact, and convex; if $h$ is also Fr\'{e}chet differentiable at $a$, then $\partial\,h(a)=\{h^{\prime}(a)\}$.
\end{itemize}

\gap

Our primary examples of Euclidean Jordan algebras are  $\Rn$, $\Sn$, and $\Hn$. In  the algebra $\Rn$ (with componentwise multiplication as Jordan product and usual inner product), spectral sets/functions (also  called symmetric sets/functions)
are precisely those that are invariant under the action of permutation matrices.
In this algebra, any two elements operator commute; strong operator commutativity requires simultaneous (permutation) rearrangement with decreasing components.
For example, in $\R^2$, the elements $(1,0)$ and $(0,1)$ operator commute, but not strongly.
In the algebras $\Sn$ (of all $n\times n$ real symmetric matrices) and $\Hn$ (of all $n\times n$ complex Hermitian matrices), the Jordan product and the inner product are given, respectively, by
$$X\circ Y:=\frac{XY+YX}{2}\quad\mbox{and}\quad \langle X,Y\rangle:=\tr(XY).$$
In $\Sn$ (in $\Hn$) spectral sets are those that are invariant under linear transformations of the form $X\mapsto UXU^*$, where
$U$ is an orthogonal (respectively, unitary) matrix.
Also, two matrices $X$ and $Y$ in $\Sn$ (in $\Hn$) operator commute if and only if
$XY=YX$, or equivalently, there exists an  orthogonal (respectively, unitary) matrix $U$ such that $X=UD_1U^*$ and $Y=UD_2U^*$, where $D_1$ and $D_2$ are diagonal matrices consisting, respectively, of eigenvalues of $X$ and $Y$. If the diagonal vectors of $D_1$ and $D_2$ have decreasing components, then $X$ and $Y$ strongly operator commute.

\gap

We now state our {\it geometric commutation principle}:

\begin{theorem}\label{geometric version}
{\it Suppose $E$ is a spectral set in $\V$ and $a\in E$. Then, $a$ strongly operator commutes with every element in the normal cone $N_{E}(a)$. 
In particular, $a$ strongly operator commutes with every element in the normal cone $N_{[a]}(a)$, where $[a]$ is the $\lambda$-orbit of $a$.} 
\end{theorem}

We prove this result as a simple consequence of (what we call) Fan-Theobald-von Neumann inequality (\ref{ftvn}) together with its equality case, see Theorem \ref{ftvn theorem} below. 
When $E$ is also convex, one may prove this as a consequence of a result on subgradients of convex spectral functions such as Theorem 5.5 in \cite{bauschke et al}. 

Based on Theorem \ref{geometric version}, we derive our commutation principles for optimization problems:

\begin{theorem}\label{main theorem}
{\it 
 Suppose $E$ is a spectral set in $\V$, $\Phi:\V\rightarrow \R$ is a spectral function, and $h:\V\rightarrow \R$.
\begin{itemize}
\item [(i)] If $a$ is an optimizer of the problem $\underset{E}{\max}\, (h+\Phi)$ and $h$ has a $E$-subgradient at $a$, then  $a$ strongly operator commutes with every element in 
$\partial_{E}\,h(a)$.
 \item [(ii)] If $E$ and $h$ are convex and $a$ is an optimizer of the problem $\underset{E}{\min}\, h$, then $a$ strongly operator commutes with the negative of some element in $\partial\,h(a)$.  
\end{itemize}
}
\end{theorem}

Specializing $h$ in the above result to a linear function leads to   an interesting consequence for variational inequality problems. To elaborate, consider a function $G:\V\rightarrow \R$ and  a set $E\subseteq \V$. Then, the 
{\it variational inequality problem} \cite{facchinei-pang}, VI$(G,E)$, is to find an element $a\in E$ such that 
$$\langle G(a),x-a\rangle \geq 0\,\,\mbox{for all}\,\,x\in E.$$ 
When $E$ is a closed convex cone, this becomes a {\it cone complementarity problem}. We now state a simple consequence of Theorem \ref{main theorem}.

\gap

\begin{corollary}\label{corollary}
{\it Suppose $E$ is a spectral set in $\V$, $\Phi$ is a spectral function, and $h:\V\rightarrow R$ is convex and Fr\'{e}chet differentiable. Let $c\in \V$ and $G:\V\rightarrow \R$. Then, the following statements hold: 
\begin{itemize}
\item [(i)] If $a$ is an optimizer of $\underset{E}{\max}\,( h+\Phi)$, then $a$ strongly operator commutes with $h^{\prime}(a)$.
\item [(ii)] If $h(x):=\langle c,x\rangle$ on $\V$ and $a$ is an optimizer of $\underset{E}{\max}\,( h+\Phi)$, then $a$ strongly operator commutes with $c$.
Moreover, the maximum value is $\langle \lambda(c),\lambda(a)\rangle+\Phi(a)$.
\item [(iii)] If $h(x):=\langle c,x\rangle$ on $\V$ and $a$ is an optimizer of $\underset{E}{\min}\, (h+\Phi)$, then $a$ strongly operator commutes with $-c$. Moreover, the minimum value is 
$\langle\, \wl(c),\lambda(a)\rangle+\Phi(a)$, where $\wl(c):=-\lambda(-c)$.
\item [(iv)] If $a$ solves VI$(G,E)$, then $a$ strongly operator commutes with $-G(a)$.
\end{itemize}
}
\end{corollary}

In our proofs, we employ  standard ideas/results from  convex analysis 
\cite{rockafellar} and the following key result from  Euclidean Jordan algebras  \cite{lim et al, baes, gowda-tao}:

\begin{theorem} \label{ftvn theorem}
{\it 
For all $x,y\in \V$,
\begin{equation} \label{ftvn}
\langle x,y\rangle\leq \langle \lambda(x),\lambda(y)\rangle.
\end{equation}
Equality holds in (\ref{ftvn})  if and only if $x$ and $y$ strongly operator commute. 
}
\end{theorem}

While our results are stated in the setting of 
Euclidean Jordan algebras for simplicity and ease of proofs, it is possible to describe these in a general setting/system. 
This general  system is formulated by turning (\ref{ftvn})  into an axiom and defining the concept of commutativity 
via the equality in (\ref{ftvn}). The precise formulation is as follows.  
A {\it Fan-Theobald-von Neumann system} (FTvN system, for short) \cite{gowda-ftvn},
is a triple $(\V,\W,\lambda)$, where $\V$ and $\W$ are  real inner product spaces and
$\lambda:\V\rightarrow \W$ is a norm preserving map satisfying the
property
\begin{equation} \label{intro ftvn}
\max\Big \{\langle c,x\rangle:\,x\in [u]\Big \}=\langle \lc,\lu\rangle\quad (\forall\,c,u\in \V),
\end{equation}
with $[u]:=\{x\in \V: \lambda(x)=\lambda(u)\}$.
This property is a combination of an inequality and a condition for equality. The inequality
$\langle x,y\rangle \leq \langle \lx,\ly\rangle$ (that comes from (\ref{intro ftvn})) is referred to as 
the {\it Fan-Theobald-von Neumann inequality} and the equality
$$\langle x,y\rangle = \langle \lx,\ly\rangle$$
defines  {\it commutativity} of $x$ and $y$ in this system. Spectral sets in this system are defined as sets of the form $E=\lambda^{-1}(Q)$ for some $Q\subseteq W$; spectral functions are of the form $\Phi=\phi\circ \lambda$ for some $\phi:W\rightarrow \R$.
Examples of such systems include \cite{gowda-ftvn}:
\begin{itemize}
\item [$\bullet$] The triple $(\V,\Rn,\lambda)$, where   $\V$ is a Euclidean Jordan algebra of rank $n$ carrying the
trace inner product with $\lambda:\V\rightarrow \Rn$ denoting the eigenvalue map. Commutativity in this FTvN system reduces to strong operator commutativity in the algebra $\V$. 
\item [$\bullet$] The triple $(\V,\Rn,\lambda)$, where $\V$ is a finite dimensional real vector space and $p$ is a real
homogeneous polynomial of degree $n$ that is  hyperbolic with respect to a vector $e\in \V$, complete and isometric,
with $\lambda(x)$ denoting the vector of
roots of the univariate polynomial $t\rightarrow p(te-x)$ written in the decreasing order \cite{bauschke et al}. 
\item [$\bullet$] The triple $(\V,\W,\gamma)$, where
$(\V,\G,\gamma)$ is a normal decomposition system (in particular, an Eaton triple)
and $\W:=\mbox{span}(\gamma(\V))$ \cite{lewis}.
\end{itemize}

Based on the property (\ref{intro ftvn}), one can show -- see Remark 1 below -- that an analog of Theorem \ref{geometric version} holds in any FTvN system. Consequently, all of our stated results can be extended to FTvN systems.   
\section{Preliminaries}
Throughout, we let  $(\V, \circ,\langle\cdot,\cdot\rangle)$ denote a Euclidean Jordan algebra of rank
$n$ with unit element $e$ \cite{faraut-koranyi, gowda-sznajder-tao}. Additionally, we assume that the inner product is the trace inner product, that is, $\langle x,y\rangle=\tr(x\circ y)$, where `tr' denotes the trace of an element (which is the  sum of its eigenvalues). In this setting, every Jordan frame in $\V$ is orthonormal and the eigenvalue map $\lambda:\V\rightarrow \Rn$ is an isometry. \\
It is well known  that any Euclidean Jordan algebra is a direct product/sum
of simple Euclidean Jordan algebras and every simple Euclidean Jordan algebra is isomorphic to one of five algebras,
three of which are the algebras of $n\times n$ real/complex/quaternion Hermitian matrices. The other two are: the algebra of $3\times 3$ octonion Hermitian matrices and the Jordan spin algebra.

It is known \cite{jeong-gowda} that {\it when $\V$ is simple, spectral sets are precisely those that are invariant under automorphisms of $\V$ (which are invertible linear transformations from $\V$ to $\V$ that preserve Jordan products).
}
For an element $a\in \V$, we abbreviate the spectral decomposition 
$a=a_1e_1+a_2e_2+\cdots+a_ne_n$ as $a=q*{\cal E}$, where $q=(a_1,a_2,\ldots, a_n)\in \Rn$ and ${\cal E}:=(e_1,e_2,\ldots, e_n)$ is an ordered Jordan frame. Note that (by rearranging the entries of $q$ and ${\cal E}$, if necessary), we can always write $a=\lambda(a)*{\cal E}$ for some ${\cal E}$.\\

We recall the following equivalent formulations of commutativity. 
\begin{proposition} (\cite{faraut-koranyi}, Lemma X.2.2; \cite{gowda-ftvn}, Prop. 2.6) \label{prop}
{\it Let $a,b\in \V$.
\begin{itemize}
\item [$(i)$] $a$ and $b$ operator commute if and only if the linear operators $L_a$ and $L_b$ commute, where $L_a(x):=a\circ x$, etc.
\item [$(ii)$] The following are equivalent:
\begin{itemize}
\item [$\bullet$] $a$ and $b$ strongly operator commute.
\item [$\bullet$] $\langle a,b\rangle=\langle \lambda(a),\lambda(b)\rangle.$
\item [$\bullet$] $\lambda(a+b)=\lambda(a)+\lambda(b).$
\end{itemize}
\end{itemize}
}
\end{proposition}

Using the fact that in a simple algebra, every Jordan frame can be mapped on to any another Jordan frame by an automorphism (\cite{faraut-koranyi}, Theorem IV.2.5), it is easily seen that: {\it If $a$ and $b$ operator commute in a simple algebra, then for some automorphism $A$ of $\V$, $A(a)$ and $b$ strongly operator commute. } A similar conclusion holds in $\Rn$ as well. 
\section{Proofs}     
\noindent{\bf Proof of Theorem \ref{geometric version}}:
Let $E$ be a (nonempty) spectral set in $\V$ and $a\in E$. Then, $E=\lambda^{-1}(Q)$ for some $Q\subseteq \Rn$. As $a\in E$,
$$[a]=\{x\in \V:\lambda(x)=\lambda(a)\}\subseteq E.$$ 
Since $N_{E}(a)\subseteq N_{[a]}(a)$, it is enough to show that $a$ strongly operator commutes with every element in $N_{[a]}(a)$. (This will also prove the second part of the theorem.)
Let $d\in N_{[a]}(a)$ so that $\langle d,x-a\rangle \leq 0$ for all $x\in [a]$. Rewriting this, we see 
\begin{equation} \label{linear inequality on the orbit}
\langle x,d\rangle \leq \langle a,d\rangle\,\,\mbox{for all}\,\,x\in [a].
\end{equation}
Now, writing the spectral decomposition of $d$ as $d=\lambda(d)*{\cal E}$ for some ordered Jordan frame ${\cal E}$, we define $x:=\lambda(a)*{\cal E}$. Then, $\lambda(x)=\lambda(a)$, $x\in [a]$, and (because every Jordan frame is orthonormal)
$\langle x,d\rangle=\langle \lambda(x),\lambda(d)\rangle=\langle \lambda(a),\lambda(d)\rangle.$
Hence, from (\ref{linear inequality on the orbit}),
$$\langle \lambda(a),\lambda(d)\rangle\leq \langle a,d\rangle.$$
From Theorem \ref{ftvn theorem}, we see the equality
 $\langle a,d\rangle=\langle \lambda(a),\lambda(d)\rangle$ and the
 strong operator commutativity of $a$ and $d$. 
$\hfill$ $\qed$

\gap

\noindent{\bf Remark 1.} In the above proof, the  part  beyond (\ref{linear inequality on the orbit})
essentially says that 
$\max \{\langle d,x\rangle:\,x\in [a] \}=\langle \lambda(d),\lambda(a)\rangle,$ which is our defining property of a FTvN system. This shows that an analog of Theorem \ref{geometric version} holds in any FTvN system.

\gap

\noindent{\bf Proof of Theorem \ref{main theorem}}:
$(i)$ Suppose $a$ solves  the problem $\underset{E}{\max}\,( h+\Phi)$, where  $E=\lambda^{-1}(Q)$ for some $Q\subseteq \Rn$ and $\Phi=\phi\circ \lambda$ for some  function $\phi:\Rn\rightarrow \R$.  Then, $a\in E$ and 
$$h(a)+\Phi(a)\geq h(x)+\Phi(x)\,\,\mbox{for all}\,\,x\in E.$$ Now,
$[a]=\{x\in \V:\lambda(x)=\lambda(a)\}\subseteq E$ and so 
\begin{equation} \label{ready to cancel}
h(a)+\Phi(a)\geq h(x)+\Phi(x)\,\,\mbox{for all}\,\,x\in [a].
\end{equation}
Since $\Phi(a)=\phi(\lambda(a))=\phi(\lambda(x))=\Phi(x)$ for all $x\in [a]$, the above expression simplifies to 
\begin{equation}\label{h inequality on the orbit}
h(a)\geq h(x)\,\,\mbox{for all}\,\,x\in [a].
\end{equation}
Now, take any $d\in \partial_{E}\,h(a)$. Then, 
$$h(x)-h(a)\geq \langle d,x-a\rangle\,\,\mbox{for all}\,\,x\in E.$$ Since $[a]\subseteq E$,  (\ref{h inequality on the orbit}) leads to 
$\langle d,x-a\rangle \leq 0$ for all $x\in [a]$, that is, $d\in N_{[a]}(a)$. By Theorem \ref{geometric version}, $a$ strongly operator commutes with $d$. 
\\
$(ii)$   Suppose $E$ is  convex (in addition to being spectral)  and $a$ solves the problem $\underset{E}{\min}\, h$, where $h$ is also convex. 
 Let $\lchi$ denote the indicator function of $E$ (i.e., it takes the value zero on $E$ and infinity outside of $E$). Then, $a$ is  a  optimizer of the (global) convex problem $\underset{\V}{\min}\, (h+\lchi)$ and so
$$0\in \partial\,(h+\lchi)(a)=\partial\,h(a)+\partial\,\lchi(a),$$
where the equality comes from the subdifferential sum formula (\cite{rockafellar}, Theorem 23.8).
Hence, there is a $c\in \partial\,h(a)$ such that $-c\in \partial\,\lchi(a).$
This $c$ will have the property that $$\langle -c,x-a\rangle \leq 0\,\,\mbox{for all}\,\,x\in E,$$ that is, $-c\in N_{E}(a)$. 
By Theorem \ref{geometric version}, $a$ strongly operator commutes with $-c$. This completes the proof.
$\hfill$ $\qed$

\gap

\noindent{\bf Remark 2.} In the proof of Item $(i)$ above, we went from 
(\ref{ready to cancel}) to (\ref{h inequality on the orbit}) by canceling the  common term $\Phi(a)$. This  type of cancellation can be carried out in certain other situations -- for example, when we consider the product $h(x)\,\Phi(x)$ with $\Phi(x)>0$ for all $x\in E$. Thus, the above proof could be modified to get  results similar to $(i)$ for  other appropriate 
combinations of $h$ and $\Phi$. 

\gap

\noindent{\bf Proof of Corollary \ref{corollary}}
$(i)$ Suppose $a$ is an optimizer of $\underset{E}{\max}\,( h+\Phi)$. As $h$ is assumed to be convex and differentiable, $h^{\prime}(a)$ is the only element in $\partial\,h(a)$. The stated assertion comes from Theorem \ref{main theorem}, Item $(i)$.\\
$(ii)$ The strong commutativity part comes from $(i)$.  Also, the maximum value is $$h(a)+\Phi(a)=\langle c,a\rangle+\Phi(a)=\langle \lambda(c),\lambda(a)\rangle+\Phi(a),$$
 where the second equality comes from Theorem \ref{ftvn theorem}. \\
$(iii)$ When $h(x)=\langle c,x\rangle$ for all $x$, and $a$ solves $\underset{E}{\min}\,( h+\Phi)$, we consider the problem $\underset{E}{\max}\,( -h-\Phi)$ and apply $(ii)$ by observing that $-\Phi$ is a spectral function. Also, the minimum value is $$-\Big (\langle -c,a\rangle-\Phi(a)\Big )=-\langle \lambda(-c),\lambda(a)\rangle+\Phi(a)=\langle \wl(c),\lambda(a)\rangle+\Phi(a).$$
$(iv)$ Suppose $a$ solves VI$(G,E)$ so that $\langle G(a),x-a\rangle\geq 0$ for all $x\in E$. Then, $a$ solves the problem $\underset{E}{\min}\, h$, where $h(x):=\langle G(a),x\rangle$ for all $x\in \V$. 
By Item $(iii)$, $a$ and $-G(a)$ strongly operator commute.
$\hfill$ $\qed$

\gap

\noindent{\bf Remark 3.} We note that strong operator commutativity of $a$ and $b$ implies the operator commutativity of $a$ and $\pm b$. Hence, Items $(ii)-(iv)$ in Corollary \ref{corollary} improve known operator commutativity relations (\cite{ramirez-seeger-sossa}, Theorem 2 and Proposition 8) for linear $h$ and variational inequalities. 
We also note that  this Corollary is similar to
Theorem 1.3 in \cite{gowda-jeong}, which is applicable to {\it simple} Euclidean Jordan algebras.\\

We now provide some illustrative examples. 

\gap
 
\noindent{\bf Example 1.}
This  example shows that in Theorem \ref{ramirez theorem},  differentibility alone is not enough to give strong operator commutativity. In the Euclidean Jordan algebra $\R^2$ spectral sets are just permutation invariant sets. So the set $E=\{(1,0),\,(0,1)\}$ is spectral. For the function
$h(x,y):=\frac{1}{2}x^2-x+x(y^2+y)$, we have $h(1,0)=-\frac{1}{2}$ and $h(0,1)=0$.
Also, $h^{\prime}(x,y)=(x-1+y^2+y,2xy+x)$. So, $h^{\prime}(1,0)=(0,1)$ and $h^{\prime}(0,1)=(1,0)$. We note that the elements $(1,0)$ and $(0,1)$ operator
commute in $\R^2$, but not strongly. Thus, if  $a$ denotes either a minimizer or a maximizer of $h$ on $E$, then $a$ and $h^\prime (a)$ do not strongly operator commute.

\gap

{\bf Example 2.} 
Consider  two $n\times n$ complex Hermitian matrices $C$ and $A$ with eigenvalues
$c_1\geq c_2\geq \cdots \geq c_n$ and $a_1\geq a_2\geq \cdots\geq a_n$. In the algebra $\Hn$, consider the spectral set
$$E:=\{UAU^*: \,U\in {\cal C}^{n\times n} \,\,\mbox{is unitary}\}.$$
As this is also compact, the linear function $\langle C,X\rangle$ attains its maximum on this set at some matrix $D$ in $E$. By Corollary 1, Item $(ii)$,
 $C$ and $D$ strongly operator commute  
and $$\underset{X\in E}{\max}\,\langle C,X\rangle =\langle C,D\rangle=\langle \lambda(C),\lambda(D)\rangle=\langle \lambda(C),\lambda(A)\rangle=\sum_{i=1}^{n}c_ia_i.$$
Thus we get the 
  classical result of Fan, namely,
$$\max \Big \{\mbox{tr}(CUAU^*):\,U\in {\cal C}^{n\times n} \,\,\mbox{is unitary} \Big\}= \sum_{i=1}^{n}c_ia_i.$$

\gap

{\bf Example 3.} In $\V$, an element $c$ is said to be an {\it idempotent} if $c^2=c$. It is known that zero and one are the only possible  eigenvalues of such an element.
If $c$ has exactly $k$ nonzero eigenvalues (namely, ones), then we say that $c$ has {\it rank} $k$.   Every idempotent of rank $k$ is of the form $e_1+e_2+\cdots+e_k$ for some Jordan frame $\{e_1,e_2,\ldots, e_n\}$.  
 Now, consider the set of all idempotents of rank $k$, where $1\leq k\leq n$.
This set is a spectral set in $\V$; it is also known to be compact. 
Now, for any $c\in \V$, we maximize $\langle c,x\rangle$ over this spectral set. By Corollary 1, the maximum is attained at some $a$ which strongly operator commutes with $c$. So, this maximum = $\langle c,a\rangle=\langle \lambda(c),\lambda(a)\rangle=\lambda_1(c)+\lambda_2(c)+\cdots+\lambda_k(c)$ since $\lambda(a)=(1,1,\ldots, 1,0,0,\ldots,0)$. Thus, {\it for any $c\in \V$, the sum of the largest $k$ eigenvalues equals the maximum of $\langle c,x\rangle$ over the set of all idempotents of rank $k$.} This is a well-known {\it variational principle}, see \cite{baes}. We remark that Theorem 1 falls short of justifying this principle.
For a broader result in the setting of certain hyperbolic systems, see 
\cite{bauschke et al}, Corollary 5.6.

\gap

{\bf Example 4.} In $\V$, let $K$ be a  closed convex cone  that is also a spectral set. For example, $K=\lambda^{-1}(Q)$, where $Q$ is a permutation invariant closed convex cone in $\Rn$ \cite{jeong-gowda}. For a function $f:\V\rightarrow \V$, consider the cone complementarity problem, CP$(f,K)$, which is to find
$x\in \V$ such that 
$$x\in K,\,y:=f(x)\in K^*,\,\mbox{and}\,\langle x,y\rangle=0,$$
where $K^*$ denotes the dual of $K$ in $\V$. We specialize Corollary 1, Item $(iv)$ to get: If $a$ solves CP$(f,K)$, then $a$ strongly operator commutes with $-f(a)$. This means that 
$$a\in K,\,\, b:=f(a)\in K^*,\,\,\mbox{and}\,\,0=\langle a,b\rangle=\langle \lambda(a),\wl(b)\rangle,$$
where $\wl(b)=-\lambda(-b)$. 


\end{document}